\documentclass[leqno,12pt]{article}
\usepackage{graphicx}
\usepackage{epsfig}
\usepackage{amsfonts,amsmath}
\usepackage{latexsym}
\usepackage{amssymb}
	\usepackage{theorem}
\usepackage{color}

\newcommand{\pf}{\hfill$\bullet$}
\usepackage[inoutnumbered,ruled]{algorithm2e}             % Algorithme
\newcommand{\demi}{\frac{1}{2}}

\def \Sup{\displaystyle\sup}
\def \Max{\displaystyle\max}

\def \Lim{\displaystyle\lim}

\newlength{\oldparindent}
\makeatletter
\newcommand*{\rom}[2]{\expandafter\@slowromancap\romannumeral #1@}

\title{Optimal Contract with Moral Hazard for Public Private Partnerships\thanks{We thank Guillaume Carlier, Ivar Ekeland, Dylan Possama\"i, Nizar Touzi and St\'ephane Villeneuve for interesting and helpful discussions and comments. }}
\author{Ishak Hajjej \thanks{\small University Tunis Manar, ENIT, Lamsin}, Caroline Hillairet \thanks{\small ENSAE ParisTech,  Crest. The author acknowledges funding from the research programs Chaire Risques Financiers of Fondation du Risque and  Investissements d'Avenir (ANR-11-IDEX-0003/Labex Ecodec/ANR-11-LABX-0047). },
 Mohamed Mnif\thanks{\small University Tunis Manar, ENIT, Lamsin } and   Monique Pontier \thanks{\small Institut de Math\'ematiques de Toulouse }
}
\date {9 Janvier 2017}

\newtheorem{definition}{Definition }[section]
\newtheorem{proposition}[definition]
{Proposition }
\newtheorem{lemme}[definition]%
{Lemma }
\newtheorem{theoreme}[definition]%
{Theorem }
{Corollary }
\newtheorem{remarque}[definition]%
{Remark }
\newtheorem{hypothese}[definition]%
{Assumption}

\newtheorem{ex}[definition]%
{Example}
\newcommand{\un}{{\mathbf{1}}} %indicatrice

\newcommand{\pr}{\mathbb{P}}

\newcommand{\ff}{\mathbb{F}}

\newcommand{\R}{\mathbb{R}}

\newcommand{\cA}{\mathcal{A}}

\def \proof{{\noindent \bf Proof: }}

\numberwithin{equation}{section}

\begin{document}

%\\
\maketitle

%\addtolength{\baselineskip}{0.3cm}

\begin{abstract}
Public-Private Partnership (PPP) is a contract between a public entity and a consortium, in which the public outsources the construction and the maintenance of an equipment (hospital, university, prison...). One drawback of this contract is that the public may not be  able to observe the effort of the consortium but only its impact on the social welfare of the project.
We aim to characterize the optimal contract for a PPP in this setting of asymmetric information between the two parties. This leads to a stochastic control under partial information and it is also related to
principal-agent problems with moral hazard. Considering a wider set of information for the public and using martingale arguments in the spirit of Sannikov \cite{San}, the optimization problem can be reduced to a standard stochastic control problem, that is solved numerically.
 We then  prove that for the optimal contract, the effort of the consortium is explicitly characterized.
In particular, it is shown that the optimal rent is not a linear function of the effort, contrary to some models of  the economic literature on PPP contracts.

\end{abstract}

{\itshape Keywords :}  Public Private Partnership, stochastic control under partial information, HJB equation, Moral Hazard.

%%%%%%%%%%%%%%%%%%%%%%%%%
\section{Introduction}
\addcontentsline{toc}{section}{Introduction}

A Public Private Parternship contract is defined by the split between private and public tasks concerning a public services, namely:
 the design of the project, the construction (building), the financing  and the maintenance (operate). DBFO means that all the four tasks are supported
by the private partner. The goal of PPP contracts is to transfer the risk to the consortium, to provide a better value for money in the use
of public funds. In France,  by the law of 2008 endorsing the order of 17th June 2004, PPP contract can not be used except if it is expressly
justified with regarded to at least one of the following criteria: emergency, complexity, economic efficiency...and actually the conclusion is that
almost all  projects are in emergency...

The relevance of outsourcing an investment in order to reduce the debt of a public entity has been  studied in Espinosa et al. \cite{CaroMo2}.
% problem of outsourcing the debt: an economic point of view, and taking into account
% the constraints that a country faces when issuing a new amount of debt, is it
% optimal for this country to finance a public project via a private investment?
Here we do not focus on the cost of the construction  but on the maintenance aspect of the PPP contract.
Hillairet and Pontier propose in \cite{CaroMo}  a study
on PPP and their relevance, assuming the eventuality of a default of the counterparty. In their model,  as in other economic papers such as Iossa et al. \cite{IMP}, the rent is assumed to be a linear rule of the effort of the consortium: although this modelisation leads to tractable computations, it seems very "ad hoc" and  economically questionable. The present work does not assume any a priori form for the rent, and  in our numerical example, it is shown that the optimal rent is actually not a linear rule.

%$\bullet$ a part on Information, Moral Hazard, principal agent literature. \\
This paper focuses on the informational asymmetry issue in  PPP contracts. Indeed, public and private partners obviously do not share the same information
for negotiation, management and follow-up of the contract.
Auriol-Picard \cite{Auriol} prove that Build-Operate-Transfer (BOT) contracts (a variant of PPP contracts) may be relevant for the public in case of better  information of the private partner, provided a large enough
number of concession candidates. But  for example  in France only three consortium  are able to support a PPP contract (Bouygues, Vinci, Eiffage).
The support mission of PPP (MAPPP in French, for Mission d'Appui aux PPP), responsible for evaluating the projects in view of legitimate the use
of a PPP contract, aims also to reduce the information asymmetry between public entity and consortium. However,  as pointed out by the  General Inspectorate of Finance in December 2012,
the multiple roles of the mission put it "de facto" in a potential situation of conflict of interest.

Besides,
the public may not be able to observe the effort of the consortium, but only its impact on the social welfare of the project.
Thus characterizing { an} optimal PPP contract in this setting of asymmetric information between both partners is related to  principal-agent problems with
moral hazard.  As shown in {book of Cvitanic et al. } \cite{bookCvitanic}, a general theory can be used to solve these problems, by means of forward-backward stochastic differential
equations.
This work is inspired by the literature on dynamic contracting using recursive methods, and in particular the
seminal paper  of Sannikov \cite{San} (2008).
In Biais et al. \cite{Biais}, the agent is risk-neutral  and his efforts, unobservable by the principal,  reduce the likelihood of large
 (but relatively infrequent) losses of the size of a project: more precisely, the  losses
occur according to a Poisson process whose intensity is controlled by the agent.
Pag\`es and Possama\"i \cite{Pages} propose an optimal contracting between competitive investors and an impatient bank monitoring a pool of long-term
loans subject to Markov contagion. The unobservable bank  monitoring decision affects the default intensity of an entity of the pool.
Optimal contracting in a  Brownian setting with risk-averse agent and principal has also been studied recently in Cvitanic et al. \cite{Touzi}, by identifying a family of admissible contracts
for which the optimal agent's action is explicitly characterized, and leading to a tractable case for CARA (exponential) utility functions.

In this paper, due to the long maturity of PPP,  we consider a perpetual contract between a public entity and a consortium. The consortium supports the initial cost of the project as well as the maintenance costs. The effort that the consortium  does  to improve the social value of the project is not observable by the public.  Thus the rent the public pays to the consortium, to compensate him for his efforts and for the operational costs,  is determined  on the basis of  the public information, that is according to the social value of the project. This is related to principal/agent problem with moral hazard and our approach relies on stochastic control under partial information, as in Bensoussan \cite{Alain}.  We consider a Stackelberg leadership model: the public (the principal) is the leader by offering a contract (characterized by the rent), while the consortium (the agent) gives a best response (characterized by the effort). The aim of this paper is to characterize such optimal contracts.
To overcome the difficulty that the control process of the consortium (the effort) is not observable by the public,
we restrict the family of admissible contracts to a set of contracts that lead to a tractable characterization of the consortium effort. This could be economically interpreted by the fact that others contracts,  for which the public does not know what incentives they will provide to the consortium effort, will likely not be offered. Moreover, we theoretically prove that the optimal contract is indeed of this form. Finally we characterize optimal contracts and provide numerical solutions.

This paper is organized as follows. Section \ref{sec2} presents the problem, Section
\ref{sec3} provides the solution of this optimal control via Hamilton-Jacobi-Belman equation.
 Section \ref{sec4} concludes with numerical illustrations
based on the Howard algorithm.

\section{Public Private Partnership's optimal contracts}
\label{sec2}

Throughout the paper,  $(\Omega,\mathbb{F}=(\mathcal{F}_t)_{t\in [0, T]}  , \mathbb{P})$  is  a filtered probability space, with $\mathbb{F}$ a Brownian filtration generated
 by  a standard Brownian  motion $W$.

\subsection{Effort and rent}

The {\bf operational cost} $(C_s)_{s \geq 0}$ of the project, supported by the consortium (and not observed by the public),   is a non-negative $\mathbb{F}$-adapted process
\begin{equation}
\label{cout}
%C_t=C_0+\int_0^t kds+\int_0^t \sigma dB_s, \sigma>0.
C_t=C_0+ k t+\sigma W_t,
\end{equation}
%{\color{blue} j'ai ecrit directement avec  $k$ et $\sigma$ constantes}
where
\begin{itemize}
 \item $C_0>0$ is the initial cost of the project, taking into account  the construction of the infrastructure.
\item  $C_t$ is the cumulative cost of the project over the period $[0,t]$, taking into account both  the cost of the construction and the cost of  the infrastructure maintenance.
\item $k>0$ and $\sigma>0$ are respectively  the drift and the volatility  of the operational cost of  the infrastructure maintenance.
\end{itemize}

\begin{remarque}
\label{coutpositif}
The cost process $(C_s)$ is not necessarily non-negative
 for all $s.$ Nevertheless,  as it is proved in Appendix \ref{app}, a sufficient condition
 to get the cost non-negative on time interval $[0,T]$
with at least probability $0.95$ is  \, $C_0/(\sigma\sqrt T)\geq 1.96.$
\end{remarque}

The consortium supports the operational cost and chooses the  effort he does to improve his service for the project : the effort is a
 non-negative $\mathbb{F}$-adapted process $(A_s)_{s \geq 0}$, it improves  the social value  of the project.
%and reduces the maintenance cost.
The social welfare, defined
as the social value of the project plus the  operational cost, is a $\mathbb{F}$-adapted process $(X_s)_{s \geq 0}$ given by
\begin{equation}
\label{eq:X}
 X_t:= X_0 + \int_0^t ( \varphi(A_s)ds + dC_s)= X_0 + \int_0^t   (\varphi(A_s)+k)ds+\sigma dW_s
 \end{equation}
where $X_0$ is the initial value of the project (i.e. of the construction, it may be a function of $C_0$) and  $\varphi$ is specified hereafter.

 The public observes  the social value  $X$ of the project, but he does not observe directly the effort of the consortium. Thus his
 information is conveyed by the filtration $\mathbb{F}^X$ generated by the social value process $X$.
The public  chooses the rent he will  pay to the consortium to compensate him for his efforts and the operational costs that he supports;
 the  rent is a  non-negative $\mathbb{F}^X$-adapted process $(R_s)_{s \geq 0}$.

 Thus  we are looking for optimal control processes $(R,A)$
with $R$ adapted to the filtration generated by the observation $X$ but itself is depending on the control
 process $A$. Remark that in our model the effort $A$ only affects the drift and not the volatility of
  the social welfare $X$ (the case of an impact both on the drift and the volatility will be done in a future work). We develop here a strong approach,
 %(instead of a weak approach as in Cvitanic et al. \cite{Touzi}, in which the volatility of the state price process is also affected by the control)
in the context of   stochastic control under partial observation as in Bensoussan \cite{Alain} Section 2.3.

% This is related to a problem of moral hazard. But, according to Bensoussan \cite{Alain} Section 2.3, we have to deal with some difficulty: we are looking for optimal control $(R,A)$ adapted to the filtration generated by the observation $X,$ but itself is depending on the control. Thus, this difficulty is relaxed in defining the controls as adapted to the filtration $\F^W$, actually generated by the observation $X$ under null control (meaning $A=0$, cf.\cite{Alain}  Section 2.3.2, (2.3.10) and (2.3.11)), so the context is the one of stochastic control under partial observation.

\subsection{A Stackelberg leadership model}
We define the respective optimization problems for the consortium and the public.
Due to the long maturity of PPP contract (up to 30-50 years), we assume that the contract is perpetual. The public and the consortium have the same time preference parameter $\delta >0$.
Let us first  define the functions involved in the formulation of the optimization problems:

\begin{hypothese}\label{Hypfct} \mbox{}\\
\begin{itemize}
\item $U$ is the utility function of the consortium, strictly concave
 strictly increasing and satisfying
$U(0)=0$ and Inada's conditions $U'(\infty)=0, U'(0)=\infty$.
\item   $\varphi$ models the impact of the consortium's efforts on the social value, $\varphi$
 is strictly concave increasing satisfying  $\varphi(0)=0$,  $\varphi'(\infty)=1$ (so $\varphi(x)\geq x$), $\varphi'(0)<\infty$.
\item  $h$ is the cost of the effort for the consortium; $h$ is convex,
$h(0)=0,$ $h'(0)>0,$ (thus $h$ is increasing) and  $h'(\infty)=\infty.$
\item   $h\circ  \psi^{-1}$ is convex where $\psi:=\frac{1}{2}(\frac{h'}{\varphi'} )^2.$
%\item  $c>0$ is the minimum  payment to the consortium, as a compensation of the investment in the construction.
\end{itemize}
 Finally, the public does not want to pay a rent over a given  amount $\bar{r}.$
 \end{hypothese}

\begin{remarque} The  function $\frac{h'}{\varphi'}$ is increasing positive, and
 $\psi(e)=\frac{1}{2}(\frac{h'}{\varphi'}(e) )^2\geq  \frac{1}{2}(\frac{h'(0)}{\varphi'(0)} )^2>0.$
\end{remarque}
%{

\noindent We define different sets of admissible contracts, depending on the information flow:
\begin{eqnarray*}
\mathcal{A}&:= &\{ (R_s,A_s)_{s\geq 0 }   \quad  \mathbb{F}\mbox{-adapted, }   A_s\geq 0, ds\otimes d\mathbb{P}\,  a.e. \mbox{ such that}
\\
& &U(R_s-k)-h(A_s)\geq 0~ ds\otimes d\mathbb{P} \, a.e. \, \mbox{ and  }   k\leq R_s\leq \bar r,~\}.
\\
\mathcal{A}^X&:=&\mathcal{A} \cap  \{ (R_s,A_s)_{s\geq 0 }, \mbox{ such that  }  \, \quad R \mbox{ is } \mathbb{F}^X\mbox{-adapted } ~\}.
\end{eqnarray*}
Those admissibility conditions ensure that  entering into the contract provides a non-negative value for the consortium.
Remark that  $(R,A)$ in $\mathcal A$  implies the following integrability properties
\begin{equation}
\label{integrable}
 e^{-\delta s }(\varphi(A_s) -R_s+k)^{-} \in L^1(\R^+\times \Omega),
 ~e^{-\delta s }(U(R_s-k)-h(A_s)) \mbox{ and }\frac{h'}{\varphi'}(A_s)\in L^2_{loc}(\R^+\times \Omega).
\end{equation}

\noindent We consider a Stackelberg leadership model: the public is the leader by offering a contract (characterized by the rent process $R$). The
 consortium gives a best response in terms of the effort process $(A_s)_{s\geq 0 }$.\\

\noindent {\bf Objective function and continuation value for the consortium and for the public:}\\
The consortium aims to optimize the expectation of his aggregate
utility of the rent minus the drift of the operational cost, minus the cost of his effort
\begin{equation}\label{objectiveC}
A \rightarrow J_0^C(R,A) =\mathbb{E} \left ( \int_0^\infty e^{-\delta s} \left( U(R_s-k) - h(A_s)\right) ds \right) .
\end{equation}
The public anticipates the consortium's best response to propose the optimal contract and  aims to optimize the expectation of the social welfare minus the rent paid to the consortium
\begin{eqnarray}\label{objectiveP}
 R \rightarrow J_0^P(R,A)&=&  \mathbb{E} \left ( \int_0^\infty e^{-\delta s} \left( dX_s-R_s ds\right)  \right)\nonumber\\
&=&\mathbb{E}\left(\int_0^\infty e^{-\delta s}( dC_s+\varphi(A_s)ds) - R_s ds)\right)\nonumber\\
 &=&\mathbb{E}\left(\int_0^\infty e^{-\delta s}( \varphi(A_s) - R_s +k) ds \right),
  \end{eqnarray}
 the last equality being a consequence of the dynamics of the social welfare
$X$ and the fact that $\mathbb{E} \left(\int_0^\infty e^{-\delta t} \sigma dW_t \right)=0$.
According to \eqref{integrable},  the integrals of both  objective
 functions   equation (\ref{objectiveC}) and equation (\ref{objectiveP}) are well defined.
From a dynamic point of view, the  objective function  at time $t$ for the consortium  is $\mathbb{P}$-a.s.
\begin{equation}
\label{maxprive}
  J_t^C(R,A):=    \mathbb{E} \left ( \int_t^\infty e^{-\delta (s-t)} \left( U(R_s-k) - h(A_s)\right) ds |\mathcal{F}_t \right)
\end{equation}
while the    objective function  at time $t$ for the public   is $\mathbb{P}$-a.s.
\begin{equation}
\label{maxpublic}
J^P_t(R,A):=    \mathbb{E}\left(\int_t^\infty e^{-\delta ( s-t)}( \varphi(A_s) - R_s +k) ds |\mathcal{F}^X_t \right).
\end{equation}
%ess\,sup_{A }
The consortium chooses the effort $A$ and the public chooses the rent $R$ in the  set of admissible contracts $\mathcal{A}^X$,
 such as to optimize their respective objective  functions,
leading to the corresponding  continuation  value process denoted respectively $V^C$ and $V^P$.\\

\noindent More precisely, an effort process $A$ is incentive compatible with respect to a given  rent $R$ if it optimizes the consortium's expected utility
 (defined in equation (\ref{objectiveC}))  given $R$. The problem of the public is to find $R$ (the contract)  that
  optimizes his expected discounted profit (defined in equation (\ref{objectiveP})), given the corresponding incentive compatible effort $A$.

 Compared to a classic optimization control problem, the difficulty of our formulation  is that the public does not observe the control $A$
of the consortium, but he observes only its impact on the social value $X$ which is the state process of the optimization control problem.
The state process $X$ appears  in an implicit way in the formulation of the optimization problem of the consortium, through the rent process
$R$, control of the public. Thus there is no explicit link between the two controls $A$ and $R$, the only indirect link involves the state
 process $X$. The trick to overcome this difficulty is to reformulate the optimization problems in terms of the   consortium continuation value process $V^C$.

\subsection{Incentive compatible contract}
To encourage the consortium to follow the recommended effort, the public proposes an incentive compatible contract. This subsection characterizes the incentive compatible contracts for $(R,A)$ in  $\mathcal{A}$ the largest set of admissibility.
As  in  Sannikov \cite{San} or  Cvitanic et al. \cite{Touzi}
in a weak formulation setting,  the following Proposition
 \ref{OptPriv} characterizes the dynamics of consortium continuation value process. It is coherent with the result of \cite{Touzi}, proved using BSDE's technics,  in a more general framework.

%\textcolor{red}{J'ai modifi\'e cette proposition
\begin{proposition}
\label{OptPriv}
If the contract $(R,A) \in \mathcal{A}$ is incentive compatible, with $A$ taking value in $]0,\bar a[$ with $\bar a := h^{-1} \circ U (\bar r-k)$, then
the dynamics of the consortium objective function is
\begin{equation*}
dJ_t^C(R,A)=\delta J_t^C(R,A) dt-( U(R_t-k) - h(A_t))dt+Y_tdW_t, \quad {\mbox a.s.}
\end{equation*}
where
 \begin{equation}
\label{Ye}
Y_t=\sigma\frac{h'(A_t)}{\varphi'(A_t)}>0.
\end{equation}
Therefore,    $A_t=(\frac{h'}{\varphi'})^{-1}\left(Y_t\sigma^{-1}\right)$,  denoted as $A^*(Y_t)$,
realizes the optimal value  in (\ref{maxprive}) for the consortium.
%the dynamics of the  consortium objective function {(in case of incentive compatible controls)}, is
%$$dJ_t^C(R,A)=\delta J_t^C(R,A) dt-( U(R_t-k) - h(A_t))dt + \sigma\frac{h'(A_t)}{\varphi'(A_t)} dW_t, \quad {\mbox a.s.}$$
If $R$ is the optimal rent for the public, then the corresponding incentive compatible effort takes value in $]0,\bar a[$ and the {\bf incentive compatible}  dynamic of the  consortium continuation value process is
 \begin{equation}\label{dWbis}
dV^C_t=\delta V_t^C  dt-( U(R_t-k) - h(A_t))dt + \sigma\frac{h'(A_t)}{\varphi'(A_t)} dW_t.
\end{equation}

\end{proposition}

\begin{remarque}\label{A*}
Incentive compatible contracts imply that the effort
$A$ is necessarily defined as an $\mathbb{F}^{V^C}=(\sigma(V^C_s,s\leq t)_{t \geq 0}$ adapted process
since it has to satisfy $A_t=(\frac{h'}{\varphi'})^{-1}\left(Y_t\sigma^{-1}\right)$ where $Y_t^2=\frac{d}{dt}\langle V^C\rangle_t.$
\end{remarque}

\proof  For any admissible pair $(R,A) \in \cA,$ let us define the process
 $$M^C_t(R,A):=e^{-\delta t}J_t^C(R,A)+\int_0^t e^{-\delta s}
\left( U(R_s-k) - h(A_s) \right) ds, \, \, \, t \geq 0 . $$
The pair $(R,A) \in \cA,$   $ U(R_s-k) - h(A_s) $ takes its values in the bounded interval
$[0,U(\bar r-k)].$ Therefore $M^C$ is an $\ff$-martingale, uniformly integrable, as
 an   $\ff$-conditional expectation of a bounded random
variable:  there exists a
 $\ff$-predictable process  $Y$ such that $e^{-\delta s}Y_s\in L^2(\Omega\times\R^+)$ and for any $t$:
\begin{equation}
\label{eq2.1}
e^{-\delta t}J_t^C(R,A)+\int_0^t e^{-\delta s} ( U(R_s-k) - h(A_s)) ds= J_0^C(R,A)+\int_0^te^{-\delta s}Y_sdW_s.
\end{equation}
The boundedness of $ U(R_s-k) - h(A_s) $ implies that the process $e^{-\delta t}J_t^C(R,A)$ is uniformly
bounded by the integral $\int_t^\infty e^{-\delta s}ds$ which goes to $0$ when $t$ goes to infinity. Thus
$t$ going to infinity in (\ref{eq2.1})  leads to the consortium's objective value
$$
\int_0^\infty e^{-\delta s} ( U(R_s-k) - h(A_s)) ds-\int_0^\infty e^{-\delta s}Y_sdW_s= J_0^C(R,A).
$$
Using  Definition (\ref{eq:X})
\begin{equation}
\label{eq2.2}
J_0^C(R,A)=\int_0^{\infty} e^{-\delta s} ( U(R_s-k)+Y_s\sigma^{-1}(\varphi(A_s)+k) - h(A_s)) ds-\int_0^{\infty}e^{-\delta s}\sigma^{-1}Y_sdX_s.
\end{equation}
 The public observes the social value $(X_t)_{t\geq 0}$ but he could not make difference between the effort $(A_t)_{t\geq 0}$ and the Brownian motion $(W_t)_{t\geq 0}$. The consortium knows that the incentive contract proposed by the public does not optimize the integral
 $\int_0^{\infty}e^{-\delta s}\sigma^{-1}Y_sdX_s$. In order to motivate the consortium, the public proposes a contract such that the corresponding optimal effort  $(A^*_t)_{t\geq 0}$  maximizes
 the concave function  $a \rightarrow Y_t \sigma^{-1}\varphi(a) -h(a),$  for all $t$, $dt \otimes d \mathbb{P}$ almost everywhere.\\
 \noindent It remains to prove that, for the "optimal" contract, the optimum on $[0,\bar a]$ of this function   is not attained on the bounds of the interval (where $\bar a := h^{-1} \circ U (\bar r-k)$ is the upper bound of the effort). More precisely, we prove that
if the incentive compatible effort is equal to $0$ or to $\bar a$, then the public could propose a better contract (that is a rent) that will  increase his value function. \\
Let $t\geq 0$ and $\eta >0$ be fixed and consider the stochastic set
$D_{t,\eta}:=\{\omega:~A^*_s(\omega)=0 \mbox{ for  }  s \in [t,t+\eta] \}$.
 From the definition of the public continuation value process,    $R^*_.=k$ on $D_{t,\eta}.$
 Therefore, on this set, $V_.^C$ is a constant process. Since the dynamics of the consortium
 continuation value process follows
 $
dV^C_s=\delta V_s^C  ds + Y_s dW_s, \mbox{ on }D_{t,\eta},
$
  the uniqueness of the It\^o decomposition implies $V_.^C=Y_.=0$ on $D_{t,\eta}$ $dt\otimes d\mathbb P$ a.e.
{   On the other hand,  on   $D_{t,\eta},$ $A=0,R=k$ implies $V_s^P=0.$
However,  since the public is the leader, he could propose a rent to the consortium satisfying $R_s-k=U^{-1}(h(A_s))$. The concavity of the functions $\varphi,$ $h^{-1},$ $U,$ yields the function
$g:x\to \varphi[h^{-1} (U(x))] -x$ is concave. Moreover it satisfies $g(0)=0,~g'(0)=+\infty$ and
going  to $-\infty$ when $x\to\infty$, thus
$\sup_{x\geq 0}(\varphi[h^{-1}(U(x)])-x)>0$ and  $V_s^P>0$.
This is a contradiction and $\mathbb P(D_{t,\eta})=0$.} This shows that the incentive compatible  effort for the "optimal" contract  satisfies $A^*_t>0$ $dt\otimes d \mathbb P$. Similarly $A^*_t< \bar a$ $dt\otimes d\mathbb P$.
%Therefore the maximum of $ a \rightarrow Y_t \sigma^{-1}\varphi(a) -h(a)$  is attained on $]0,\bar a[$  and the  relation  $Y_t=\sigma\frac{h'(A_t)}{\varphi'(A_t)}>0$ is stated.
 \pf
\\

In the following, {\bf  all admissible contracts  $(R,A)$ are assumed incentive compatible}, and thus are denoted $(R,A^*(Y))$. With a slight abuse of notation we keep the same notations $\mathcal{A}$ and $\mathcal{A}^X$. In the next section, we will first solve the problem under the set of controls $\mathcal{A}$, that is with no restriction of measurability on the rent process.
As suggested by Remark \ref{A*},  we will then check that the optimal controls $(R^*,A^*)$ over $\mathcal{A}$ are functions of the consortium continuation value $V^C$.
 Thus equation (\ref{dWbis}) modeling  the dynamics of the consortium value function is a Markovian diffusion, the solution of which being defined up to its explosion time $\tau.$ We will prove that $\tau=+\infty $ a.s.
%Thus  $(R,A)\in \mathcal{A}^M$ are both defined on $[ 0,\tau ].$

 %According to a measurable section theorem (see \cite{castaing2}),
 The    public  value at time $0$  over the class ${\cal A}$ is written as  follows
\begin{equation}\label{**}
v(x):=\sup_{(R,A=A^*(Y))\in {\cal A}}J_0^P(R,A)=\sup_{(R,A=A^*(Y))\in {\cal A}} E_x [\int_0^\infty  e^{-\delta s} (\varphi(A_s)-R_s+k)ds]
\end{equation}
where  $E_x$ is the  conditional expectation with respect to the event $\{ V_0^C= x \}$.
%{\color{red} j'ai remis tous les indices $M$ dans les $\cal A$ qui avaient tous disparus!? On resout bien d'abord sur $A^M$ et ensuite on verifie que c'est en fait aussi le meilleur controle sur $A^X$}.
 The fact that the effort $A=A^*(Y)$ is the best response of the consortium, for a given rent $R$, follows from  the incentive compatible dynamic (\ref{dWbis}) of the state process $V^C$.
\\

The next section characterizes,  through a HJB equation, the function $v$ that realizes the optimum for the public over the class of ${\cal A}$.
%Nevertheless, this does not change the consortium/public optimal values  in the general setting
But, actually we will prove that  the optimal processes  over the class $\mathcal{A}$
are in fact $\mathbb{F}^X$-adapted and
the optimal value function for the public over the class $\mathcal{A}^X$ is
indeed $v( V^C_0) $ (see Proposition \ref{optimumFXadapted}).
Thus we solve ultimately the original optimization problem over the class $\mathcal{A}^X$.

\section{Optimal controls and value functions for the public and the consortium}
\label{sec3}
We first solve the optimization problem (\ref{**}) over the set of controls $\mathcal{A}$, that is done using the dynamic programming principle.
Subsection \ref{formalHJB} gives a formal derivation of the Hamilton Jacobi Bellman equation
\begin{eqnarray}\label{hjbra}
\sup_{(r,a)\in{\cal C}}\left[-\delta w(x)+{\cal L}^{r,a}w(x) -r+\varphi(a)+k\right]=0,
\end{eqnarray}
where the second order differential operator ${\cal L}^{r,a}$ is defined by
$${\cal L}^{r,a}w(x):= w'(x)(\delta x-U(r-k) + h(a) )+\demi w''(x)\left(\sigma\frac{h'}{\varphi'}(a)\right)^2$$
and the control space ${\cal C}$ is defined by
\begin{equation}
\label{defCalC}
{\cal C}:=\{(r,a) \in [k,\bar{r}]\times \mathbb{R}^+,~h(a)\leq U(r-k)\}.
\end{equation}
Subsection \ref{proofHJB} relates the HJB equation to the optimization problem (\ref{**}), and characterizes the optimal controls and value functions  over $\mathcal{A}$.
We then prove in Subsection \ref{A=AM} that the optimal controls over the larger set $\mathcal{A}$ are actually in $\mathcal{A}^X$.
\subsection{Formal derivation of the  HJB equation}
\label{formalHJB}
We assume that the public value function $v$ (defined in equation \eqref{**}) is  of class  $C^2$.
Standard methods (e.g. El Karoui \cite{NEK} or Pham \cite{HP}) are used
to provide the HJB equation that  is satisfied by
the public value function $v$. We give nevertheless some details for the sake of completeness.
\\
Let $\eta>0,$ and $(R,A^*(Y))\in \mathcal{A}$, % and the stopping time $\tau^C=\inf\{t \geq 0,~V^C_t\leq 0\}$.
the dynamic programming principle yields
$$v(x)\geq \mathbb{E}_x\left[\int_0^{\eta}  e^{-\delta s}(\varphi(A_s^*(Y))-R_s+k)ds+
e^{-\delta \eta} v(V^C_{\eta})\right].$$
%where $E_x$ is the  conditional expectation with respect to the event $\{ V^C_0= x \}$.
It\^o's formula applied to the process $(e^{-\delta s}v(J^C_s(R,A^*(Y))))_{s\geq 0}$ and taking expectation yield
$$0\geq \mathbb{E}_x[\int_0^{\eta} e^{-\delta s}(-\delta v(x) +{\cal L}^{R_s,A^*_s(Y)}v(J^C_s)+\varphi(A_s^*(Y))-R_s+k)ds].$$
Dividing by $\eta,$ letting $\eta$ goes to $0$, then using the continuity of the integrands and the mean value theorem we obtain
$$-\delta v(x)+v'(x)(\delta x-U(r-k) + h(a) )+
\demi v"(x) ( \sigma\frac{h'(a)}{\varphi'(a)} )^2 -r+\varphi(a)+k\leq 0.$$
Since this holds for any control $(r,a)\in {\cal C},$ we obtain the inequality
\begin{equation}
\label{ineq}
\sup_{(r,a)\in{\cal C}}\left[-\delta v(x)+{\cal L}^{r,a}v(x) -r+\varphi(a)+k\right]\leq 0.
\end{equation}
On the other hand, suppose that $(R^*,A^*(Y))$ is an optimal control for the class $\mathcal{A}$. Then the dynamic
programming principle yields
$$v(x)=\mathbb{E}_x\left[\int_0^{\eta} e^{-\delta s}({\cal L}^{R^*_s,A_s^*(Y)}v(J^C_s)+\varphi(A_s^*(Y))-R^*_s+k)ds\right]$$
where $V^C_s$ is the consortium value function at time $s$, given by  the optimal control
$(R^*,A^*(Y))$. By similar arguments, dividing by $\eta,$ and sending $\eta$ to $0$,
one has at time $t=0$
$$-\delta v(x)+{\cal L}^{R^*_0,A_0^*(Y)}v(x)+\varphi(A_0^*(Y))-R^*_0+k=0$$
which combined with (\ref{ineq}) yields (\ref{hjbra}).

\noindent
{ A useful result to solve this HJB equation is the following.}
\begin{lemme}
\label{vBounded}
The function $v$ defined in \eqref{**} is a non-negative bounded function on $\mathbb{R}^+$.
\end{lemme}
\proof
 For all $(R,A) \in \cal A$ and for any $t \geq 0$,  $J^C_t(R,A)\geq 0 $ thus $V^C_t\geq 0.$ Therefore the function $v$ is defined  on $\mathbb R^+$.\\
For any $(r,a)\in{\cal C},$ $ h(a) \leq U(r-k)$  and due to the concavity of $\varphi$, $\varphi(a)-(r-k) \leq \varphi'(0)h^{-1}( U(r-k))-(r-k)$.
Since the function
$x\to \varphi'(0)h^{-1} (U(x)) -x$ is concave,
going from $0$ at $x=0$  to $-\infty$ when $x\to\infty$, it admits the maximum
\begin{equation}
\label{supv}
\sup_{x\geq 0}(\varphi'(0)h^{-1}(U(x))-x)
\end{equation}
 which bounds $v$ from above on $\mathbb{R}^+$.
Besides,   the constant control $R_t=k; A_t=0,\,\,dt\otimes d\mathbb{P}, \, a.e., $ is admissible and incentive compatible (the best consortium's response to a minimum rent is zero effort). Using this control implies that $\forall x\geqslant0, v(x)\geq 0.$
\pf
\\

\subsection{Verification theorem}\label{proofHJB}
In this subsection the optimal control problem (\ref{**}) is characterized as the solution of the HJB equation
 (\ref{hjbra}) via a verification theorem (cf. for instance   El Karoui \cite{NEK}, Pham \cite{HP}, Krylov
 \cite{Krylov} or Fleming-Rishel \cite{FR}). The following proposition provides the structure of the optimal control.

\begin{proposition}
\label{existMax}
Suppose % {\bf$(H_{U,\varphi,h,\psi})$}
Assumption \ref{Hypfct}, then there exists a unique admissible optimal  pair $(r^*,a^*)$ which realizes the maximum in HJB equation (\ref{hjbra}).
This  optimal pair  is defined  by $r^*(x)=(k+(U')^{-1}(\frac{-1}{w'(x)}) )\un_{w'(x)<0} \wedge \bar{r}$
and $a^*(x)=\arg\max (a\to  w'(x) h(a)+w"(x) \sigma^2\psi(a) +\varphi(a)) $
in the  compact interval $\left[0,h^{-1}\circ U\left((U')^{-1}(\frac{-1}{w'(x)})\un_{w'(x)<0}\wedge (\bar{r}-k)\right)\right].$ Moreover, the function $r^*$ is continuous on  $\R^+$ and $a^*$ is Borel measurable on $\left[0,h^{-1}\circ U\left((U')^{-1}(\frac{-1}{w'(x)})\un_{w'(x)<0}\wedge (\bar{r}-k)\right)\right]$.
%as soon as  $w$ is $C^3$-class: il reste encore un point a verifier pour appliquer
%le th des fct implicites...}
\end{proposition}

\begin{ex}
\label{3.3}
For instance,  the hypotheses of Proposition \ref{existMax} are satisfied for:
$\varphi(x)=x+\ln(1+x),$
$h(x)=\frac{2}{3}\sqrt{1+x}(x+4)-\frac{8}{3},$ so $h'(x)=\frac{2+x}{\sqrt{1+x}},$
$\varphi'(x)=1+\frac{1}{1+x}$, $2\psi(x)=1+x,$ $\psi^{-1}(x)=2x-1.$
The functions  $\varphi$  and $\psi$ are concave,  and $\psi^{-1}$ is convex. Finally
$h"(x)=\frac{1}{\sqrt{1+x}}-\frac{2+x}{2(1+x)^{3/2}}=
\frac{x}{2(1+x)^{3/2}}\geq 0,$ $h$ is convex, $h(0)=0$,
$h'(0)=2$, $h$ convex yields $h(x)\geq 2x.$\\
In this example, the coefficient of $w''(x)$ is $\demi\sigma^2(1+a)$, so it does not depend on $x.$
\end{ex}
\proof
(i) {\it Change of variables}:
$$P=h(a)-U(r-k)~;~Q=\demi(\frac{h'(a)}{\varphi'(a)} )^2=\psi(a).$$
Recall $\psi:x\to\demi(\frac{h'(x)}{\varphi'(x)} )^2$ is non decreasing since $h$ is convex and $\varphi$ concave.
This change of variable is a bijection
 $$ [k,\bar{r}]\times\R_+\to [-U(\bar{r}-k),0]\times\R_+~; (r,a)\mapsto (h(a)-U(r-k),\psi(a)).$$
with inverse
$$[-U(\bar{r}-k),0]\times[\psi(0),\psi(\infty)) \to [k,\bar{r}]\times\R_+ ~;~(P,Q)\mapsto ( k+U^{-1}[h\circ\psi^{-1}(Q)-P], \psi^{-1}(Q)).$$
%{\color{blue}il faut verifier une derniere fois les $-c$ ci dessous}
The constraint $P\leq 0\Leftrightarrow U(r-k)-h(a)\geq 0$ is the admissibility condition on  control $(r,a).$
\\
Thus  the HJB equation \eqref{hjbra} is equivalent to

{\small
$$\sup_{-U(\bar{r}-k) \leq P\leq 0,\psi(0)\leq Q\leq \psi(\infty)}\left[w'(x)P+w"(x)\sigma^2Q+\varphi(\psi^{-1}(Q))-
U^{-1}(h\circ\psi^{-1}(Q)-P))\right]=-\delta xw'(x)+\delta w(x).$$}

\noindent (ii) {\it
For any  $C^2$ function $f$  and under the  hypotheses of Proposition \ref{existMax}, the function
$${\cal H}:~(P,Q)\mapsto f'(x)P+f"(x)\sigma^2Q+\varphi(\psi^{-1}(Q))-U^{-1}(h\circ\psi^{-1}(Q)-P{})+k$$
is strictly concave.}\\
To prove that  ${\cal H}$ is strictly concave, we check that $D^2(-\cal{H})(P,Q)$ is a  positive-definite matrix:
 $$ \partial^2_{P,P}{\cal H}< 0,\mbox{ and }\partial^2_{Q,Q}{\cal H} < 0
 \mbox{ and } \partial^2_{P,P}{\cal H}\partial^2_{Q,Q}{\cal H}-(\partial^2_{P,Q}{\cal H})^2> 0.$$
Computing the first order derivatives:
\begin{eqnarray}
\label{der1}
&&\partial_{P}{\cal H}=(U^{-1})'(h\circ\psi^{-1}(Q)-P)+f'(x),
\\
\partial_{Q}{\cal H}&=&
\varphi'(\psi^{-1}(Q))(\psi^{-1})'(Q)-
(U^{-1})'(h\circ\psi^{-1}(Q)-P).
(h\circ\psi^{-1})'(Q)+\sigma^2 f"(x).\nonumber
\end{eqnarray}
Using the concavity of function $U$ the second  order derivative $\partial^2_{P,P}{\cal H}$ satisfies
$$
\partial^2_{P,P}{\cal H}=-(U^{-1})"(h\circ\psi^{-1}(Q)-P))<0.
$$
Other concavity arguments yield
\begin{eqnarray*}
\partial^2_{Q,Q}{\cal H}&=&
(\psi^{-1})"(Q)\varphi'\circ\psi^{-1}(Q)+((\psi^{-1})'(Q))^2\varphi"\circ(\psi^{-1})(Q)\\
&-&(h\circ\psi^{-1})"(Q)(U^{-1})'(h\circ\psi^{-1}(Q)-P)\\
&-&((h\circ\psi^{-1})'(Q))^2
(U^{-1})"(h\circ\psi^{-1}(Q)-P)<0.
\end{eqnarray*}
%puisque $\psi$ convexe montre $\psi^{-1}$ concave, que $U$ concave montre $U^{-1}$
%convexe, que $\varphi$ est concave,
%et que par hypoth\`ase $h\circ\psi^{-1}$ est convexe.
%\\
Finally
$$\partial^2_{P,Q}{\cal H}=(h\circ\psi^{-1})'(Q)(U^{-1})"(h\circ\psi^{-1}(Q)-P)$$
and since
$(U^{-1})"(h\circ\psi^{-1}(Q)-P)>0$, the Jacobian sign  is the one of
$$
-(\psi^{-1})"(Q)\varphi'\circ\psi^{-1}(Q)-((\psi^{-1})'(Q))^2\varphi"\circ(\psi^{-1})(Q)+(h\circ\psi^{-1})"{(Q)}(U^{-1})'(h\circ\psi^{-1}(Q)-P)
$$
which is positive  (using  again concavity arguments). Then
%puisque $\varphi$ et $\psi^{-1}$ sont concaves et $h\circ\psi^{-1}$ convexe.
%Alors
 $$\partial^2_{P,P}{\cal H}\partial^2_{Q,Q}{\cal H}-(\partial^2_{P,Q}{\cal H})^2>0$$
and    $\cal H$ is strictly concave.

\noindent (iii) {\it Existence of an optimal pair  $(r^*,a^*)$}.\\
In the  HJB equation one has to maximize the function
$$g: r\mapsto -w'(x)U(r-k)-r.$$
When $w'(x)\geq 0,$  this function is non-increasing and the optimum is $k.$\\
Otherwise, the function is concave and  the optimum is achieved for $$r^*(x)=(k+(U')^{-1}(\frac{-1}{w'(x)})\un_{w'(x)<0})\wedge \bar{r}.$$
Moreover $x\to r^*(x)$  is continuous even
at zero points of the function   $w'$ because $(U')^{-1}(\infty)=0$.\\
The optimal pair has to satisfy the admissibility condition\\
$$0\leq h(a^*)\leq U\left((U')^{-1}(\frac{-1}{w'(x)})\un_{w'(x)<0}\wedge (\bar{r}-k)\right).$$
 Since the function
 $$a\mapsto w'(x) h(a)+w"(x) \sigma^2\psi(a) +\varphi(a) $$
is continuous on the compact interval $\left[0,h^{-1}\circ U\left((U')^{-1}(\frac{-1}{w'(x)})\un_{w'(x)<0}\wedge (\bar{r}-k)\right)\right]$
 there exists an optimal solution $a^*$. Moreover by a selection theorem (See Appendix B, Fleming and Rishel \cite{FR}), there exists a Borel-measurable function $a^*$ from $\left[0,h^{-1}\circ U\left((U')^{-1}(\frac{-1}{w'(x)})\un_{w'(x)<0}\wedge (\bar{r}-k)\right)\right]$ into ${\mathbb R}$.
Finally, by construction $(r^*,a^*)$ takes its values in $\cal C$.
\pf
\\
\begin{proposition}\label{prop:intervalleVo}
Define  $\bar{x}:=\frac{1}{\delta}\left( U o (U')^{-1} (\frac{h'(0)}{\varphi'(0)})\right)$ and suppose Assumption  \ref{Hypfct}. Then
\\
(i) The consortium initial value $V_0^C$ is in the interval $[0, \bar{x}]$.
\\
(ii) $ v(\bar{x})=0.$
\end{proposition}
\proof
i) We assume $V_0^C > \frac{1}{\delta}(U o (U')^{-1}( \frac{h'(0)}{\varphi'(0)}))$ and we prove that it leads to a contradiction.
Let $r_0:=k+U^{-1}(\delta V^C_0)$
 and define the strict convex function $\tilde h(.):= h(.)-U'(r_0-k) \varphi(.)$. Then
%using
%$$\tilde h'(0)>U'(r-k)\Leftrightarrow U'(r-k)<\frac{h'(0)}{1+\varphi'(0)}\Leftrightarrow
%r>(U')^{-1}(\frac{h'(0)}{\varphi'(0)})+k$$
$\delta V_0^C > U o (U')^{-1} (\frac{h'(0)}{\varphi'(0)})$,  $U$ increasing and $U'$ decreasing  imply that $U'(r_0-k)\varphi'(0)<  h '(0)$.
\noindent For any control $(R,A^*(Y)) \in \mathcal{A}$, one has
by  strict concavity of the function $U$
\begin{equation}\label{inegalite11}
 U(R_s-k)
< U(r_0-k)+  U'(r_0-k)(R_s-r_0 ), \quad ds \otimes d\mathbb{P} \quad a.e.~~\mbox{ when}~ R_s\ne r_0.
\end{equation}
and by strict  convexity of $\tilde h(.):= h(.)-U'(r_0-k) \varphi(.)$,  with our choice of $r_0$, $h(0)=\varphi(0)=0, \tilde{h}'(0)>0$,
\begin{equation}\label{inegalite22}
-h(A_s^*(Y))< - U'(r_0-k)  (\varphi(A_s^*(Y))) \quad ds \otimes d\mathbb{P} \quad a.e.\mbox{ when }~ A_s> 0.
\end{equation}
For any control $(R,A^*(Y)) \in \mathcal{A}$ that is not identical to $(r_0,0)$,
\begin{eqnarray*}
J^C_0(R,A^*(Y)) &=&\mathbb{E}( \int_0^{\infty} e^{-\delta s} (  U(R_s-k)-h(A_s^*(Y)) )ds)\\
&<& \mathbb{E}( \int_0^{\infty} e^{-\delta s} ( U(r_0-k)+  U'(r_0-k)(R_s-\varphi(A_s^*(Y))-r_0 )ds)\\
&<& \frac{U(r_0-k)}{\delta} - U'(r_0-k) ( \frac{r_0-k}{\delta} + \mathbb{E}( \int_0^{\infty} e^{-\delta s} (\varphi(A_s^*(Y))-R_s+k)ds)),
%\mbox{ when} A_s> 0 ~\mbox{and} ~R_s\ne r
\end{eqnarray*}
 where the last inequality holds from \eqref{inegalite11} and  \eqref{inegalite22}.
 Since $r_0=U^{-1}(\delta V^C_0 )+k$, for $(R,A) \neq (r_0,0)$,
\begin{equation}\label{eq:inegalitestrict}
\delta J^C_0(R,A)<\delta V^C_0- U'(r_0-k)[r_0-k+\delta J^P_0(R,A)].
 \end{equation}
% $(*<)$ est vrai sauf si $A_s=0$ and  $R_s=r=k+U^{-1}(V_0^c)>k \quad ds \otimes d\mathbb{P} \quad a.e$ ce n'est pas le contr\^ole optimal car  si la réponse optimal de consortium ($A=0$) le public pr\'ef\'rera proposer une rente $R=k$\\
%Since $(r,0)$ can not be the optimal control as $r>k$, the
%optimal control $(R^*,A^*)$ satisfies the strict inequality (\ref{eq:inegalitestrict})  and using  $J_0^C(R^*,A^*)=V_0^C$,  and $V_0^P=v(V_0^C)$, we get
By taking the supremum over all admissible strategies, we get
 \begin{equation}\label{majorationVPP}
 \delta V_0^P  \leq k-r_0< 0.
 \end{equation}
On the other hand, by Lemma \ref{vBounded}, one has $V_0^P \geq 0$. This gives a contradiction to  \eqref{majorationVPP}, and $V_0^C$ is
necessarily smaller than $\frac{1}{\delta}\left( U o (U')^{-1} (\frac{h'(0)}{\varphi'(0)})\right)$. \\

ii)  We now prove that $v(\bar x)=0$.\\
The proof follows  the same lines as the previous part $i)$ with this time $V^C_0=\bar x= \frac{1}{\delta}( U o (U')^{-1}( \frac{h'(0)}{\varphi'(0)}))$.\\
We notice that the constant control $(r_0,0)$, which leads to $J_0^C(r_0,0)=\bar x$, is an admissible control and so the equality $V_0^C=\bar x$ makes sense.
%and  once again $r=k+U^{-1}(\delta V^C_0).$ Then
%using
%$$\tilde h'(0)>U'(r-k)\Leftrightarrow U'(r-k)<\frac{h'(0)}{1+\varphi'(0)}\Leftrightarrow
%r>(U')^{-1}(\frac{h'(0)}{\varphi'(0)})+k$$
Let us recall that  $r_0=k+U^{-1}(\delta V^C_0)$,  $U'(r_0-k)\varphi'(0) =  h '(0)$.
\noindent For any control $(R,A^*(Y)) \in \mathcal{A}$, one has
by  concavity of the function $U$
\begin{equation*}\label{inegalite1}
 U(R_s-k)
\leqslant  U(r_0-k)+  U'(r_0-k)(R_s-r_0 ), \quad ds \otimes d\mathbb{P} \quad a.e.
\end{equation*}
and by  convexity of $\tilde h$
\begin{equation*}\label{inegalite2}
-h(A_s^*(Y))\leqslant  - U'(r_0-k)  (\varphi(A_s^*(Y))) \quad ds \otimes d\mathbb{P} \quad a.e.
\end{equation*}
thus
\begin{equation*}
J^C_0(R,A^*(Y))\leqslant \frac{U(r_0-k)}{\delta} - U'(r_0-k) ( \frac{r_0-k}{\delta} + \mathbb{E}( \int_0^{\infty} e^{-\delta s} (\varphi(A_s^*(Y))-R_s+k)ds))
\end{equation*}
%
%\begin{eqnarray*}
%J^C_0(R,A) &=&\mathbb{E}( \int_0^{\infty} e^{-\delta s} (  U(R_s-k)-h(A_s) )ds)\\
%&\leqslant & \mathbb{E}( \int_0^{\infty} e^{-\delta s} ( U(r-k)+  U'(r-k)(R_s-\varphi(A_s)-r )ds)\\
%&\leqslant& \frac{U(r-k)}{\delta} - U'(r-k) ( \frac{r-k}{\delta} + \mathbb{E}( \int_0^{\infty} e^{-\delta s} %(\varphi(A_s)-R_s+k)ds)),
%\end{eqnarray*}
 %where the last inequality holds from \eqref{inegalite1} and  \eqref{inegalite2}.
 and $r_0=k+U^{-1}(\delta \bar x)$ implies
 $$\delta J^C_0(R,A^*(Y))\leqslant\delta \bar x- U'(r_0-k)[r_0-k+\delta J^P_0(R,A^*(Y))].$$
% \\ Or $(*<)$ appliqué à $(R^*,A^*)$ contrôle optimaux
%{\color{blue} j'ai essaye de reprendre un peu la fin de redaction pour la rendre plus claire }
%By Proposition  \ref{existMax}, there exists $ (R^*,A^*),$
%such that $J^C_0(R^*,A^*)=V^C_0=\bar x,$ and $J^P_0(R^*,A^*)=v(V^C_0)=v(\bar x),$
 %thus
 Taking the supremum over all incentive compatible admissible strategies, we get
\begin{equation}\label{majorationVP}
\delta v(V^C_0)\leq k-r_0\leq 0.
\end{equation}
Besides  $v\geqslant 0$ by Lemma \ref{vBounded}. Therefore  $v(\bar x )=0$.
\pf
\begin{lemme}\label{lemmabound} The function $v$ satisfies
\begin{eqnarray}\label{cl0}
v(0)=\frac{1}{\delta}\sup_{0\leq x\leq\bar{r}-k}(\varphi\circ h^{-1}\circ U(x)-x)
\end{eqnarray}
%$$v(0)=\frac{1}{\delta}(\varphi \circ h^{-1}\circ U(x^*)-x^* $$
%where $x^*=x_0
\end{lemme}
\proof %{ Dans toute la section 3.2, on est restreint aux controles markoviens, donc prop 3.2 est utilisable dans ce lemme}
Since every admissible control $(R,A^*(Y))\in \cal A$ must satisfy $U(R_s-k)-h(A^*(Y)_s)\geq 0,\,ds \otimes d\mathbb{P}\, a.e$, then
the equality $V^C_0=0$ holds for every control $(R^*,A^*(Y))\in \cal A$ satisfying $U(R^*_s-k)-h(A^*(Y)_s)= 0,\,ds \otimes d\mathbb{P}\, a.e$ which is equivalent to
\begin{eqnarray*}
h^{-1}\circ U(R^*_s-k)=A^*(Y)_s,\,ds \otimes d\mathbb{P}\,a.e.
\end{eqnarray*}
In this class of controls, one has
\begin{eqnarray*}
\varphi \circ h^{-1}\circ U(R^*_s-k)-R^*_s+k\leq \sup_{k\leq r \leq \bar{r}}  \left( \varphi\circ h^{-1}\circ U(r-k) - (r-k)\right) \,ds \otimes d\mathbb{P}\,a.e.
\end{eqnarray*}
This shows that
\begin{eqnarray*}
V_0^P\leq \frac{1}{\delta }\sup_{k\leq r \leq \bar{r}}  \left( \varphi\circ h^{-1}\circ U(r-k) - (r-k)\right),
\end{eqnarray*}
and the equality holds when $(R^*_s,A^*(Y)_s)=(\hat r,\hat a)\,ds \otimes d\mathbb{P}$ a.e. where $\hat r$ realizes the maximum of the function
$r\longrightarrow \left( \varphi\circ h^{-1}\circ U(r-k) - (r-k)\right)$ on $[k, \bar{r}]$ and $\hat a =h^{-1}\circ U(\hat r-k)$. We easily  check that the control $(\hat r,\hat a)$ is indeed admissible (incentive compatible).
This proves (\ref{cl0}).

\hfill\pf
\\
We now proceed to the verification Theorem (cf. \cite{HP} Th. 3.5.3., infinite horizon).
Thanks to Proposition \ref{prop:intervalleVo}  and Lemma \ref{lemmabound}, we study the HJB equation on a bounded domain $[0,\bar x]$ with the Dirichlet  boundary conditions
\begin{equation}
\label{boundHJB}
v(0)=\frac{1}{\delta}\sup_{0\leq x\leq\bar{r}-k}(\varphi\circ h^{-1}\circ U(x)-x) , \, \, \,  \, \, \,v(\bar x)= 0.
\end{equation}

\begin{theoreme}
\label{th:existsol}
Under Assumption  \ref{Hypfct}
%{\bf$(H_{U,\varphi,h,\psi})$},
the  HJB equation (\ref{hjbra}) with boundary conditions \eqref{boundHJB}
admits  a unique  solution $w$ in $C^2((0,\bar{x}))\cap C^0([0,\bar{x}])$.
Let $(a^*(x),r^*(x))$ be the argmax in the HJB equation (\ref{hjbra}). We assume that
$a^*$ is of bounded variation, then the associated stochastic differential equation
\begin{equation}\label{SDEV*}
dV^C_t=\delta V^C_tdt-(U(r^*(V^C_t)-k) - h(a^*(V^C_t)) dt+\sigma\frac{h'(a^*(V^C_t)}{\varphi'(a^*(V^C_t))} dW_t,~V^C_0=x
\end{equation}
admits a unique strong solution denoted as $V^{C}.$
We define the following controls:\\
$A^*_t:=a^*(V^C_t),R^*_t:=r^*(V^C_t)$ $dt\otimes d\pr$ almost everywhere.
Then, $(R^*,A^*)$ is the optimal control in $\mathcal{A}$. As a conclusion, the public value function satisfies $v(x)=w(x).$
\end{theoreme}

Remark that on the numerical simulations, $a^*$ is decreasing and thus is indeed of bounded variation.

\proof
\underline{Step 1}: The equation (\ref{hjbra}) admits an unique $C^2((0,\bar{x}))\cap C^0([0,\bar{x}]$ solution.\\
This follows from  Theorem 1 in \cite{Strulovici}, whose  assumptions are satisfied. Indeed the set of controls ${\cal C}=\{(r,a) \in [k,\bar{r}]\times \mathbb{R}^+,~h(a)\leq U(r-k)\}$ is a nonempty compact set. The coefficients of the HJB equation  (\ref{hjbra}) are affine functions of the variable  $x$, and continuous functions in the control $(r,a)$, for each $x$, thus standard linear growth assumptions on the coefficients of  (\ref{hjbra})  are satisfied. Under Assumption \ref{Hypfct}, the  function $\frac{h'}{\varphi'}$ is increasing positive, and
$0< \frac{h'(0)}{\varphi'(0)}\leq \frac{h'(a)}{\varphi'(a)}$.
Therefore the volatitity coefficient of the HJB equation (\ref{hjbra}) admits a uniform lower bound. In addition, the coefficient $-r+\varphi(a)+k$ is bounded on ${\cal C}$.
 Therefore, applying Theorem 1 in \cite{Strulovici},  the HJB equation (\ref{hjbra}) has a twice continously differentiable solution in $(0,\bar{x})$, and  continuous on $[0,\bar{x}]$.\\
Furthermore
 Proposition \ref{existMax} proves that for all $x$
there exists  an admissible %locally Lipschitz
 pair $(r^*(x),a^*(x))$ in $\mathcal{C}$  such that
 $$-\delta w(x)+{\cal L}^{(r^*(x),a^*(x))}w(x)- r^*(x) + \varphi(a^*(x))+ k
 =\max_{(r,a)\in{\cal C}}\left[-\delta w(x) +{\cal L}^{r,a}w(x)-r + \varphi(a)+k)\right].$$

\underline{Step 2: $w(x)\geq v(x)$}
\\
 Let $(R,A^*(Y))\in {\cal A}$ and the corresponding process $(J^C_t(R,A^*(Y)))_{t\geq 0}$ with dynamic
$$dJ^C_t(R,A^*(Y))=\delta J^C_t(R,A^*(Y))dt-(U(R_t-k) - h(A_t^*(Y))) dt+\sigma\frac{h'(A_t^*(Y))}{\varphi'(A_t^*(Y))} dW_t,~J^C_0(R,A^*(Y))=x.
$$
Applying It\^o's formula to the process $(e^{-\delta t}w(J^C_t(R,A^*(Y))))_{t \geq 0}$
{\footnotesize
\begin{eqnarray*}
e^{-\delta T}w(J^C_{T}(R,A^*(Y)))=
w(x)&+&\int_0^{T}e^{-\delta s}[-\delta w(J^C_s(R,A^*(Y)))+{\cal L}^{R,A^*(Y)}w(J^C_s)(R,A^*(Y))]ds\\
 &-&\int_0^{T}e^{-\delta s}
\sigma w'(J^C_s(R,A^*(Y)))\frac{h'(A_s^*(Y))}{\varphi'(A_s^*(Y))}dW_s
\end{eqnarray*} }
with $ w'(J^C_t(R,A^*(Y))) \frac{h'(A_t^*(Y))}{\varphi'(A_t^*(Y))}$  being a bounded process.
Taking the expectation,  for all $T$
$$E_x[e^{-\delta T}w(J^C_{T}(R,A^*(Y)))]=
w(x)+E_x[\int_0^{T}e^{-\delta s}[{\cal L}^{R_s,A_s^*(Y)}w(J^C_s(R,A^*(Y)))-\delta w(J^C_s(R,A^*(Y))) ]ds].$$
From the  HJB equation (\ref{hjbra}), ${\cal L}^{R_s,A_s^*(Y)}w(J^C_s(R,A^*(Y))) -\delta w(J^C_s(R,A^*(Y))) \leq R_s-\varphi(A_s^*(Y))-k$ thus
$$E_x[e^{-\delta T}w(J^C_{T}(R,A^*(Y)))]\leq
w(x)+E_x[\int_0^{T}e^{-\delta s}(-\varphi(A_s^*(Y))-k+R_s)ds.$$
Using boundedness from above of $-\varphi(a)-k+r$ when $(r,a)\in{\cal C}$ (cf. the proof of Lemma
\ref{vBounded})  and admissibility conditions on $(R,A^*(Y))\in {\cal A}$, one has
$e^{-\delta s}|-\varphi(A_s)-k+R_s|\leq e^{-\delta s}(\varphi(A_s)+k-R_s)^-+e^{-\delta s}M\in L^1(\R^+\times\Omega)$.
By the  dominated convergence theorem, we obtain
\begin{eqnarray*}
  \Lim_{T\longrightarrow \infty}E_x[\int_0^{T}e^{-\delta s}(-\varphi(A_s^*(Y))-k+R_s)ds]=
  E_x[\int_0^{\infty}e^{-\delta s}(-\varphi(A_s^*(Y))-k+R_s)ds].
\end{eqnarray*}
Besides, as $w$ is a continuous function  and $J^C_T$ is bounded, one has
$$\lim_{T\longrightarrow \infty} E_x[e^{-\delta T}w(J^C_T(R,A^*(Y)))]= 0\mbox{ for all }(R,A^*(Y))\in  {\cal A}.$$
Therefore \begin{eqnarray*}
0\leq w(x)+E_x[\int_0^{\infty}e^{-\delta s}(-\varphi(A_s^*(Y))-k+R_s)ds],
\end{eqnarray*}
so for any  $(R,A^*(Y))\in {\cal A}$ we get
$$w(x)\geq E_x[\int_0^{\infty}e^{-\delta s}(\varphi(A_s^*(Y))+k-R_s)ds],$$
and
$$w(x)\geq\sup_{(R,A^*(Y))  \in  {\cal A}} V^P_0(R,A^*(Y))=v(x).$$

\underline{Step 3: the SDE \eqref{SDEV*}  admits a unique strong solution}
\\
Let us consider
the SDE \eqref{SDEV*} associated to the optimal pair $(r^*(V^C_t),a^*(V^C_t))$, then $\tilde V^C_t:=e^{-\delta t}V^C_t$, $dt\otimes d \mathbb P$ a.e.
satisfies the SDE
\begin{equation}
\label{VC*}
d\tilde V^C_t=-e^{-\delta t}(U(r^*(V^C_t)-k) - h(a^*(V^C_t)) dt+\sigma e^{-\delta t}\frac{h'(a^*(V^C_t)}{\varphi'(a^*(V^C_t))} dW_t,~V^C_0=x.
\end{equation}
 The existence of a strong solution of this SDE is given by Nakao \cite{NAKAO} (Theorem p. 516),
as the drift is bounded measurable, the volatility  $\sigma\frac{h'}{\varphi'}(a^*)$ is strictly bounded from below and
 of bounded variation (since it is the case for $a^*$).  The existence of a strong solution of the SDE \eqref{SDEV*} follows.

\underline{Step 4:  $w(x)\leq v(x)$}
\\
This solution actually is the process denoted as $(V^{C}_t)_{t \geq 0}$ meaning $V^C=J^C(R^*,A^*)$ where $R^*_s:=r^*(V^C_s),~A^*:=a^*(V^C_s),\,ds\otimes d  \mathbb P\,a.e.$.
We now repeat the above arguments of Step 2:
{\footnotesize
$$e^{-\delta T\wedge\tau_n}w(V^C_{T\wedge\tau_n})=
w(x)+\int_0^{T\wedge\tau_n}e^{-\delta s}[-\delta w(V^C_s)+{\cal L}^{R^*_s,A^*_s} w(V^C_s)]ds-\int_0^{T\wedge\tau_n}e^{-\delta s}\sigma w'(V^C_s)\frac{h'(A^*_s)}{\varphi'(A^*_s)}dW_s.$$
}
Thus for $n$ and all $T$
$$E_x[e^{-\delta T\wedge\tau_n}w(V^C_{T\wedge\tau_n})]=w(x)+E_x[\int_0^{T\wedge\tau_n}e^{-\delta s}[-\delta w(V^C_s)+{\cal L}^{R^*_s,A^*_s} w(V^C_s)]ds].$$
But since $w$ satisfies the HJB equation (\ref{hjbra}) with such optimal controls we get
$$E_x[e^{-\delta T\wedge\tau_n}w(V^C_{T\wedge\tau_n})]=w(x)+E_x[\int_0^{T\wedge\tau_n}e^{-\delta s}(R^*_s-k-\varphi(A^*_s))ds]$$
$$=w(x)-E_x[\int_0^{T\wedge\tau_n}e^{-\delta s}(R^*_s-k-\varphi(A^*_s))^-ds]+E_x[\int_0^{T\wedge\tau_n}e^{-\delta s}
(R^*_s-k-\varphi(A^*_s))^+ds].$$
Taking into account the boundedness of the controls  $(r^*(V^C_t),a^*(V^C_t))_{t\geq 0}$ (cf. Proposition \ref{existMax}), the dominated convergence theorem {allows} us to get $n$ and $T$ going to infinity in $E_x[\int_0^{T\wedge\tau_n}e^{-\delta s}(R^*_s-k-\varphi(A^*_s))ds]$. Besides,
 Fatou's lemma and boundedness of $w$ allow us to
 get $n$ and $T$ going to infinity in $E_x[e^{-\delta T\wedge\tau_n}w(V^C_{T\wedge\tau_n})]$.
Therefore
$$w(x)\leq J^P_0(R^*,A^*)\leq v(x).$$
%{\color{blue}il nous semble qu'il faut developper la limite du dernier terme  de l'equation: on fait CV dominee pour la partie negative de $(R-k+\phi(A))$ et lemme de FAtou pour la partie positive. de plus, cela ressert un peu loin, en l'occurence p15 fin de la preuve de Prop3.5}\\
Conclusion:  $w$ is the public value function $v$ defined in  equation \eqref{**}
and $(a^*(V^C_t),r^*(V^C_t))$ is a Markovian optimal control.
\pf
\\

\subsection{Going back to the original set $\mathcal{A}^X$ of control processes}\label{A=AM}
In this subsection, we will prove that the process $V^{C}$ is ${\mathbb F}^X$-adapted.

\begin{proposition}
\label{optimumFXadapted}
 Under the assumptions  of Theorem \ref{th:existsol},
 the  unique strong solution to the stochastic differential equation \eqref{SDEV*} admits an infinite explosion time, and
 the filtrations $\mathbb{F}^X$,  $\mathbb{F}^{V^C}$  and $\mathbb{F}$ coincide at the optimum.
\end{proposition}

\proof
First, as $R$  and $A$ are bounded, the strong solution of the SDE \eqref{SDEV*} does not explode. Besides,
 $\mathbb{F}^{V^C}$ and $\mathbb{F}^X$  are obviously included in $\mathbb{F}$.\\
(i) We first  express $V^{C}$ as  the solution of the  SDE
\eqref{SDEV*}
%driven by $X_t$ defined in (\ref{eq:X}).
%and  to use Theorem 38  in Protter \cite{protter}(Chapter V, Section 7 p.309)  for the existence of  strong solution of EDS.
$$dV^{C}_t=
\delta V^{C}_tdt-(U(r^*(V^{C}_t)-k) - h(a^*(V^{C}_t))) dt+\sigma\frac{h'}{\varphi'}(a^*(V^{C}_t)) dW_t.
$$
Under the assumptions of Theorem \ref{th:existsol}, this SDE admits an unique strong solution, thus the filtrations generated by
$W$ and $V^{C}$ coincide since $\sigma\frac{h'}{\varphi'}(a^*(V^C_t))$ the coefficient of $dW_t$  is positive
{(cf. Corollary 1.12 of  Revuz-Yor \cite{RevYor}).}
\\
(ii) On the other hand, by definition:
$$dX_t= (\varphi(a^*(V^{C}_t))+k)dt+\sigma dW_t,$$
thus, since $\sigma>0,$
 $$dW_t=\frac{1}{\sigma}[dX_t- (\varphi(a^*(V^{C}_t))+k)dt]$$
and
$$dV^{C}_t=
\delta V^{C}_tdt-(U(r^*(V^{C}_t)-k) - h(a^*(V^{C}_t))dt+\frac{h'}{\varphi'}(a^*(V^{C}_t))\left(
- (\varphi(a^*(V^{C}_t))+k)dt +dX_t\right).
$$
Once again, under the assumptions  of Theorem \ref{th:existsol}, $V^{C}$ is a strong solution of this stochastic differential equation driven
by $X$ so this process $V^{C}$ is ${\mathbb F}^X$-adapted.
Therefore the  three filtrations $\mathbb{F}^{V^C}$, $\mathbb{F}^X$ and $\mathbb{F}$ coincide at the optimum.

\pf

%\noindent {\bf Filtrations of the information consortium/prive :}\\
            %
%
%{\color{blue}  La remarque \ref{A*} montre que $\ff=\ff^{V^C}$ des que $A$ optimal :
%pour une paire de contr\^oles optimale $(a^*,r^*)$ il vient
%\\
%$$dV_t^C(a^*,r^*)=\delta V_t^C(a^*,r^*) dt-( U(r^*) - h(a^*))dt-\sigma\frac{h'(a^*)}{1+\varphi'(a^*)}dW_t$$
%en sous entendant que dans l'EDS $(a^*,r^*)=(a^*,r^*)(V_t^C).$
%\\
%Puis apr\`es le theorem de verification, toujours \`a l'optimum,
%$V^P_t=F(V_t^C)$
%v\'erifie par la formule de It\^o
%$$dV^P_t=F'(V_t^C)dV^C_t+\demi \sigma^2\left(\frac{h'(a^*(V_t^C))}{1+\varphi'(a^*(V_t^C))}\right)^2F"(V_t^C)dt.$$
%
%
%A priori on a seulement $\ff^{V^P} \subset  \ff^X\subset\ff$, mais vu la r\'egularit\'e
%des coef de l'EDS fonctions de $V^C_t$ et le fait que $F''(V_t^C)<0$
%doit montrer qu'en fait on a aussi la coincidence entre les filtrations
%$\ff^{V^P}$ et $\ff^{V^C}.$
%}

%\newpage

\section{Numerical implementation}
\label{sec4}
The consortium continuation value  is the state parameter used in the resolution of the stochastic control  problem of Section \ref{sec3}. Proposition \ref{prop:intervalleVo} gives an upper bound for the consortium initial value $V_0^C \in [0,\bar{x}]$.

\subsection{Howard's Algorithm}

The HJB-equation \eqref{hjbra} is written as follows
\begin{equation}\label{a1}
\Sup_{(r,a)\in{\cal C}}\left[-\delta v(x)+v'(x)(\delta x-U(r-k) + h(a) )+
\demi v"(x) ( \sigma\frac{h'(a)}{\varphi'(a)} )^2 -r+\varphi(a)+k\right]= 0.
\end{equation}
Let  $\Delta$ be the finite difference step on the state coordinate and  $(x_i)_{i=1,N}$, $x_i=i\Delta$,
be the points of the grid $\Omega_\Delta.$
 The equation \eqref{a1} is discretized by replacing the first and second derivatives of $v$ with  the following approximations
\begin{displaymath}
 v'(x)\simeq \left\{ \begin{array}{ll}
\frac{v(x+\Delta)-  v(x)}{\Delta} & \textrm{ if $\delta x+h(a)-U(r-k)\geq 0$}\\
\frac{v(x)- \tilde v(x-\Delta)}{\Delta} & \textrm{ if not}
\end{array} \right.
\end{displaymath}
\begin{displaymath}
 v{''}(x) \simeq \frac{ v(x+\Delta)- 2 v(x)+ v(x-\Delta)}{\Delta^2}.
\end{displaymath}
% Taking into account the boundedness of the function $v$ and the boundary conditions in $0$, we deduce from the definition of $\tilde J$
$$
 v(0)=\frac{1}{\delta}\sup_{0\leq x\leq\bar{r}-k}\varphi\circ h^{-1}\circ U(x)-x,~v(\bar{x})=0.$$
 where $\bar{x}=\frac{1}{\delta}\left(U o (U')^{-1} (\frac{h'(0)}{\varphi'(0)}) \right)=N\Delta$.\\
 This leads to the system of $(N-1)$ equations with
 $(N-1)$ unknowns $(v^\Delta(y_i))_{i=1,...N-1}$:
  $$\max_{(r,a)\in{\cal C}}\left[A^{\Delta,(r,a)}v^\Delta(x_i)+B^{\Delta,(r,a)}\right]=0$$
where
$B^{\triangle,(r,a)}$ is given by
\begin{displaymath}
B^{\triangle,(r,a)}=
\left( \begin{array}{cccccc}
-r+\varphi(a)+k +(\frac{b^-(\Delta)}{\Delta}+\frac{a(\Delta}{\Delta^2})v(0) \\
-r+\varphi(a)+k \\
\vdots\\
 \vdots\\
-r+\varphi(a)+k\\
\end{array} \right)
\end{displaymath}
the matrix $A^{\Delta,(r,a)}$ is defined as follows:
$$[A^{\Delta,(r,a)}]_{i,i-1}=\frac{b^-(y_i)}{\Delta}+\frac{a(x_i)}{\Delta^2};~
[A^{\Delta,(r,a)}]_{i,i}=c(y_i)-\frac{|b(x_i)|}{\Delta}-2\frac{a(x_i)}{\Delta^2};
[A^{\Delta,(r,a)}]_{i,i+1}=\frac{b^+(x_i)}{\Delta}+\frac{a(x_i)}{\Delta^2};
$$
with  $b^+(x)=\Max(b(x),0)$, $b^-(x)=\Max(-b(x),0)$ and
{
\begin{eqnarray*}
 c(x)&=& -\delta
 \\
b(x)&=&( h(a)-U(r-k) ) +\delta x,
\\
a(x)&=&\demi( \sigma\frac{h'(a)}{\varphi'(a)} )^2 .
\end{eqnarray*}
}
To solve the latter equation we use  an iterative Howard algorithm  (cf. Howard \cite{howard} chapter 8).
It consists in computing two sequences $(r^n(x_i),a^n(x_i))_{i=1,...N-1}$
  and $(v^{\Delta,n}(y_i))_{i=1,...N-1}$ (starting from
$(r^1,a^1)$ chosen arbitrary):
\\
$\bullet$ step $2n-1$: to the strategy $(r^n,a^n)$ we compute $v^{\Delta,n}$  solution of  the linear system
$$A^{\Delta,(r^n,a^n)}v^{\Delta}+B^{\Delta,(r^n,a^n)}=0$$
on the grid $\Omega^\Delta.$
\\
$\bullet$ step $2n$: $v^{\Delta,n}$ is associated with a strategy
$$(r^{n+1},a^{n+1})\in\arg\max_{(r,a)\in{\cal C}} (A^{\Delta,(r,a)}v^{\Delta,n}+B^{\Delta,(r,a)}).
$$
The convergence of the Howard algorithm holds when the matrix $A^{\Delta,(r,a)}$ satisfies the
 discrete maximum principle: a sufficient condition is that  $A^{\Delta,(r,a)}$  is diagonally dominant. This  is the case since $c(x)<0$.

\subsection{Effort and rent }
We recall  $r^*(x)=(k+(U')^{-1}(\frac{-1}{v'(x)})\un_{v'(x)<0})\wedge \bar{r}$\\
and $a^*(x)=\arg\max (a\to  v'(x) h(a)+v"(x) \sigma^2\psi(a) +\varphi(a)) $\\
%$[0,h^{-1} \circ U  \left((k+(U')^{-1}(\frac{-1}{v'(x)})\un_{v'(x)<0})\wedge \bar{r}-k\right)]$ \\
in the  compact interval %$[0,  h^{-1}(U((k+(U')^{-1}(\frac{-1}{v'(x)})\un_{v'(x)<0})\wedge \bar{r}-k)].$\\
$\left[0,h^{-1}\circ U\left((U')^{-1}(\frac{-1}{v'(x)})\un_{v'(x)<0})\wedge (\bar{r}-k)\right)\right]$.

If $v'(x)\geq 0$
 $r^*(x)=k, \mbox{ and}~~ a^*(x)=\arg\max (a\to  v'(x) h(a)+v"(x) \sigma^2\psi(a) +\varphi(a)) $
in the  compact interval $[0, h^{-1}(0)]=\{0\}.$

If $v'(x)<0$
$r^*(x)=(k+(U')^{-1}(\frac{-1}{v'(x)}))\wedge \bar{r}, \mbox{ and}~~ a^*(x)=\arg\max (a\to  v'(x) h(a)+v"(x) \sigma^2\psi(a) +\varphi(a))$  in the compact interval
 $\left[0,h^{-1} \circ U  \left(((U')^{-1}(\frac{-1}{v'(x)})\un_{v'(x)<0})\wedge (\bar{r}-k)\right)\right].$

\subsection{Numerical results}

\begin{figure}[ht]
\begin{center}
\centering \includegraphics[width=9.5cm,height=11cm]{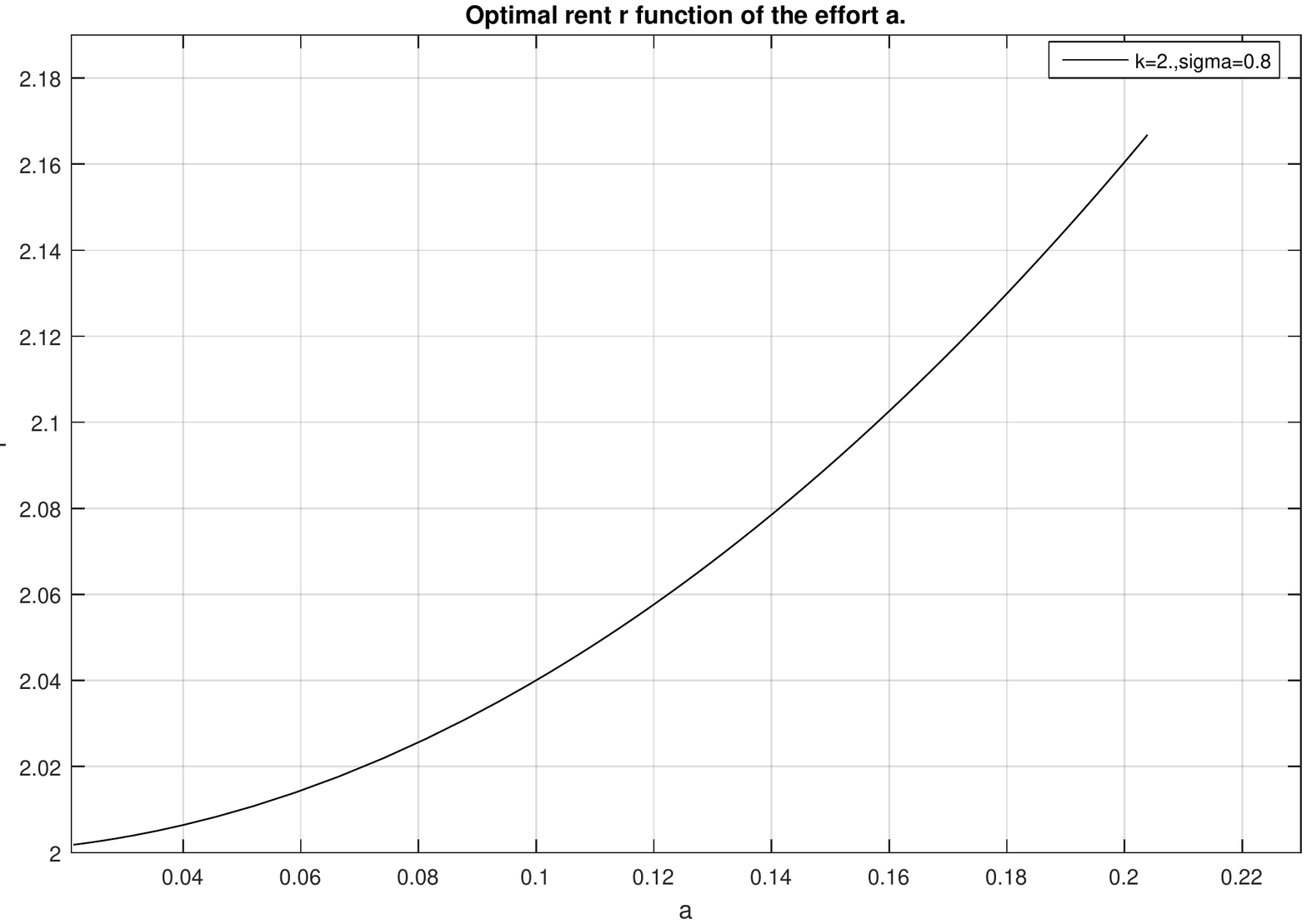}
\caption{ Optimal rent $r^*$ function of the effort  $a^*$. }
\label{RfunctionA}
\end{center}
\end{figure}

In this section we choose $U(x)=\sqrt x$  and the functions of the example \ref{3.3} :
$\varphi(x)=x+\ln(1+x),$
$h(x)=\frac{2}{3}\sqrt{1+x}(x+4)-\frac{8}{3},$ so $2\psi(x)=1+x.$\\
We choose the parameter $\delta= 0.1.$
We restrict our figures to  $x\in\left[0, \frac{1}{\delta}\left(U o (U')^{-1} (\frac{h'(0)}{\varphi'(0)}) \right)\right]$
(cf. Proposition \ref{prop:intervalleVo}), that is $[0,5]$ in our numerical example.\\
In Figure \ref{test1u}, Figure \ref{test2u} and Figure \ref{test3u}, $k=2$ and $\sigma$ varies : $\sigma=0.5, 0.8$ or $1$.
%In Figure \ref{testk1u}, Figure \ref{testk2u} and Figure \ref{testk3u}, $\sigma=2.25$ and $k$ varies : $k=1.75, 2$ or $2.3$.
Figure \ref{RfunctionA} gives the optimal rent function of the optimal effort, for $k=2$ and $\sigma=0.8$.\\
We observe that $v$ and $a^*$  are  non-increasing functions of the consortium value.

\subsubsection{Graph of the optimal rent $r$ as a function of $a$}
%{ c'est $r$ fonction de $a$ plutot qu'il faudrait}}
The more interesting observation  is that the optimal rent is an increasing convex function of
the optimal effort. This  contradicts the usual assumption of linear dependence
between rent and effort.
Besides, the {qualitative} behaviour of these optimal parameters
is the same with respect to both $\sigma$ and $k.$

%\newpage

\newpage

\subsubsection{Sensibility of the results to parameter $\sigma$ }
\begin{figure}[h!]
\begin{center}
 \includegraphics[width=9.5cm,height=9cm]{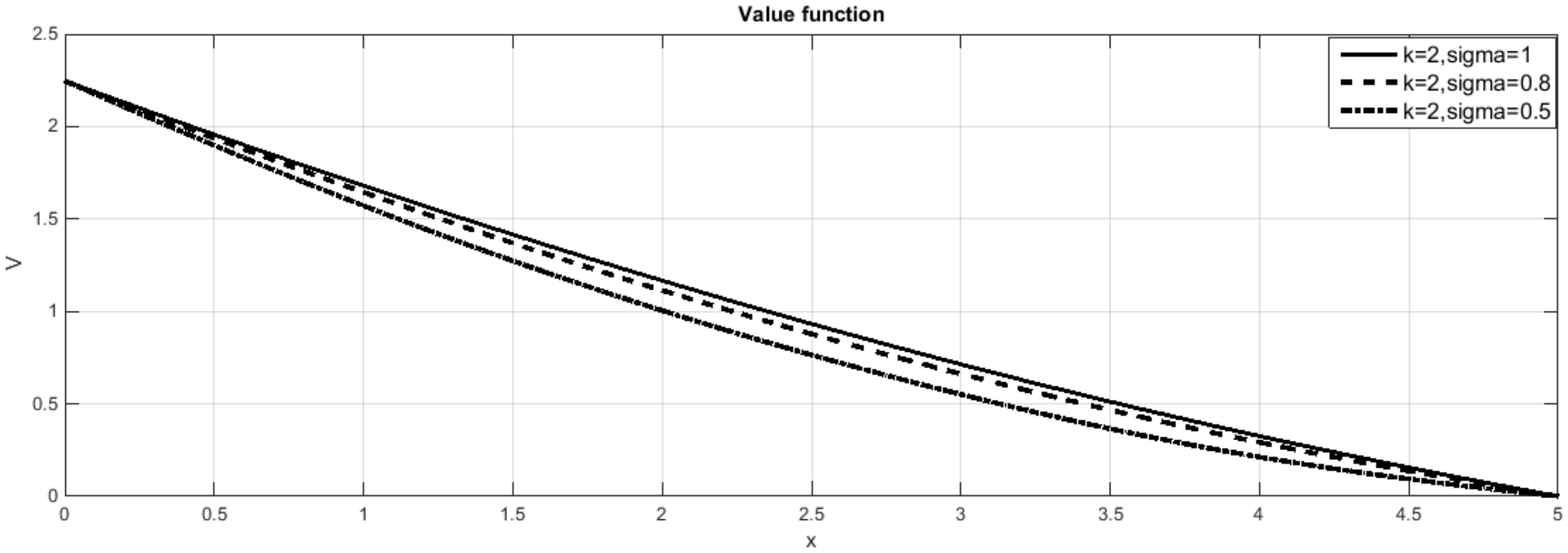}
\caption{ Value function $v$ for different $\sigma$. }
\label{test1u}
\end{center}
\end{figure}

\begin{figure}[h!]
\begin{center}
\centering \includegraphics[width=9.5cm,height=9cm]{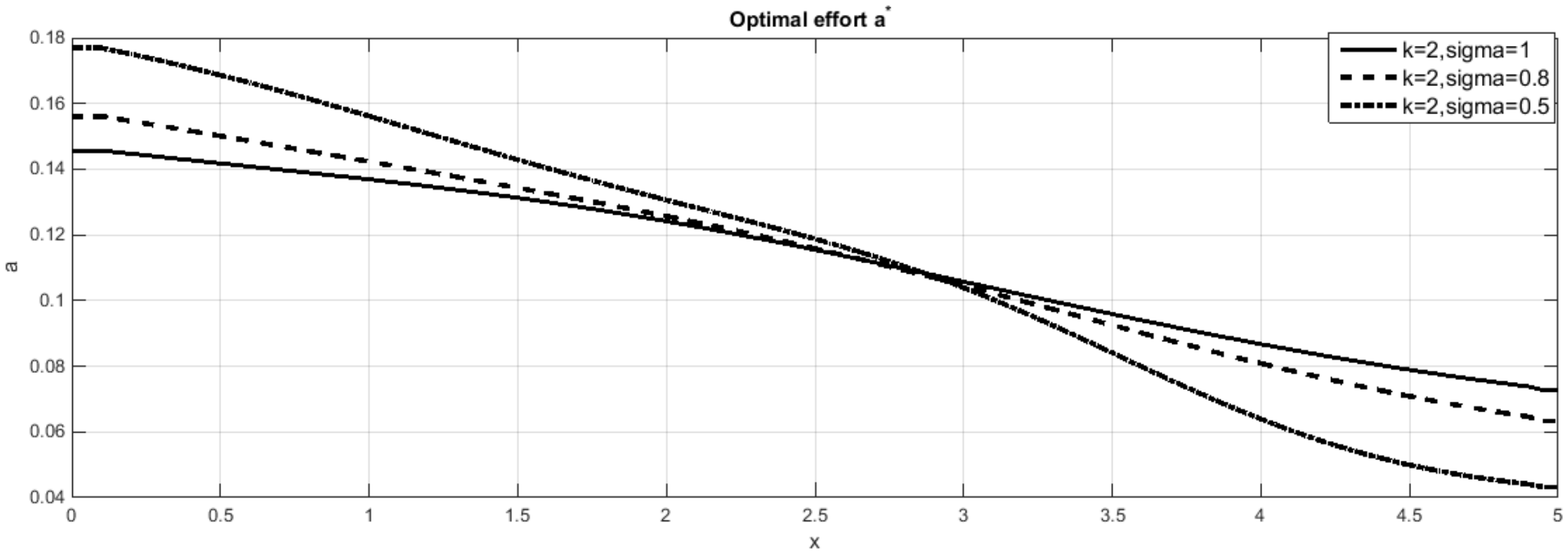}
\caption{ optimal effort $a^*$ for different $\sigma$. }
\label{test2u}
\end{center}
	\end{figure}

\begin{figure}[h!]
\begin{center}
\centering \includegraphics[width=9.5cm,height=9cm]{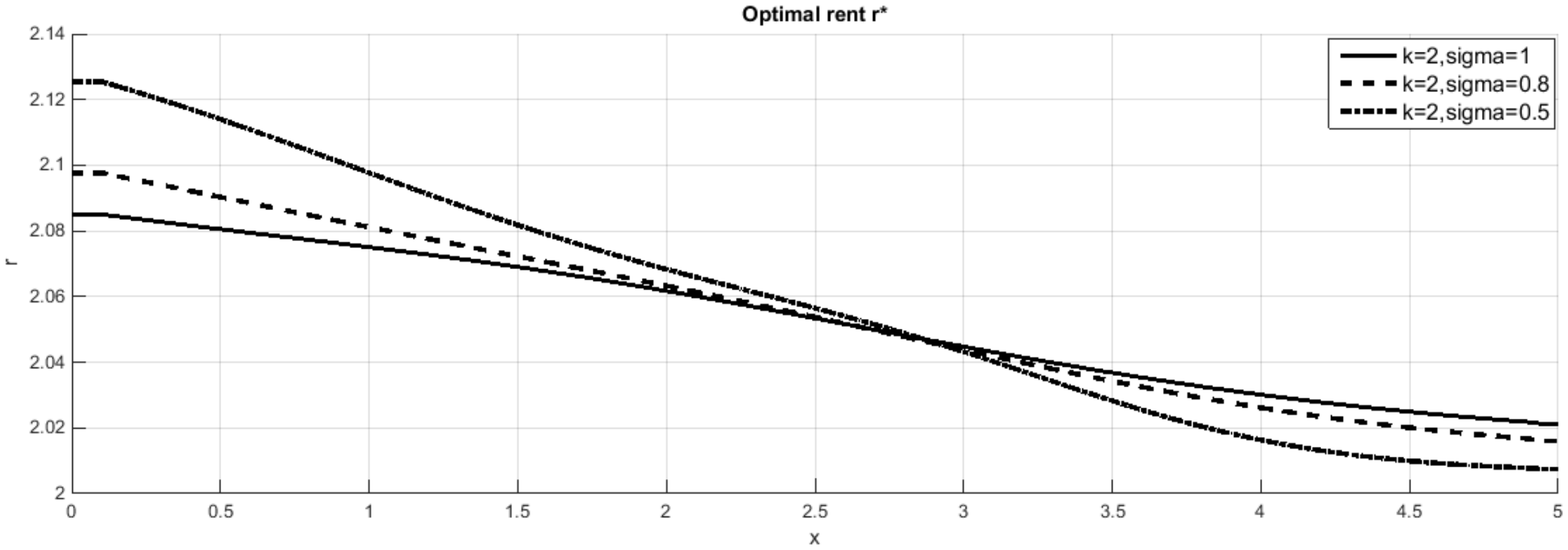}
\caption{ Optimal rent $r^*$, for different $\sigma$. }
\label{test3u}
\end{center}
\end{figure}

 According to Figure \ref{test1u}, it seems that the
optimal public value function is  increasing with respect to $\sigma$:
the  risk is supported by the consortium.
 The same behavior is observed
 for the optimal effort (Figure \ref{test2u}) and for the optimal rent (Figure \ref{test3u})
 in case of $x$ large enough. But, when $x$ is lower the optimal effort is decreasing.
 In case of a low level of the private continuation value $x$, the consortium is not ready to provide more efforts.
 This behavior is observed for any parameter $k.$
\\

%\newpage

%$$k=1,~\sigma=1,2,3,4,~c=2,~\delta= 0.1.$$
%We check that the function $v$ is positive and bounded from above.

%
%
%
%\vspace{25cm}
%
%
%\newpage

%\newpage

\subsubsection{Sensibility of the results to parameter $k$}
Actually, the true parameter control is $r-k$, so, as it could be expected, the
parameter $k$ has no impact on the behaviors of $v,$ $r,$ and $a.$
%The value function, and both optimal effort and optimal rent,
%are increasing with parameter $k.$ Indeed, we could remark that the dependence of $r^*(x)$ with respect
%to $v$ is done via $v'(x)$: $r^*(x)=k+(U')^{-1}(\frac{-1}{v'(x)})\un_{v'(x)<0}$, meaning that $r^*(x)$ is increasing with respect to $|v'(x)|$. So this behaviour is expected. As we noticed it above,
%$a^*(x)$ is also increasing with parameter $k.$

%{\color{red} le comportement des controles parait logique, mais comment interpreter celui de $v$ en fonction de $k$ ???}

%
%\begin{figure}[ht]
%\begin{center}
% \includegraphics[scale=0.35]{testk1u.pdf}
%\caption{ Value function $v$ for different $k$. }
%\label{testk1u}
%\end{center}
%\end{figure}
%
%
%
%
%\begin{figure}[ht]
%\begin{center}
%\centering \includegraphics[scale=0.30]{testk2u.pdf}
%\caption{ optimal effort $a^*$ for different $k$. }
%\label{testk2u}
%\end{center}
%\end{figure}
%
%
%
%\begin{figure}[ht]
%\begin{center}
%\centering \includegraphics[scale=0.30]{testk3u.pdf}
%\caption{ Optimal rent $r^*-k$, for different $k$. }
%\label{testk3u}
%\end{center}
%\end{figure}
%
%\vspace{35cm}
%
\newpage

\subsection*{Conclusion}
This paper provides a characterisation of optimal public private partnership contracts in a moral hazard framework,
 using martingale methods and stochastic control. A numerical example shows that, in particular, the optimal rent is
  a convex (and not linear) function of the effort. This convexity, due to the information asymmetry  between the
   consortium and the public entity, implies for the public entity a more and more costly contract to encourage
   the consortium to do more efforts. This feature should be taken into account in the models concerning PPP contracts.

%\newpage

\
%\newpage

\section{Appendix A}
\label{app}

We need to know some sufficient conditions to get the cost process $(C_s)$
non-negative on a time interval $[0,T]$. Recall the Inverse Gaussian law $(IG(\lambda,\mu))$
with density on $\R^+$
$$f(t,\mu,\lambda)=\left(\frac{\lambda}{2\pi t^3}\right)^{\demi}\exp(-\frac{\lambda}{2\mu^2 t}(t-\mu)^2).$$
The cost  process $C$ is a drifted Brownian motion and the event
$$\{\inf_{0\leq s\leq T}C_s\geq 0\}=\{\inf_{0\leq s\leq T}ks+\sigma W_s\geq -C_0\}=
\{\sup_{0\leq s\leq T}(-ks-\sigma W_s)\leq C_0\}=\{T_{C_0}> T\}$$
where $T_{C_0}$ is an hitting time. It is well known (cf. \cite{CMJ} for instance)
that the law of $T_{C_0}$ is $IG(-\frac{C_0}{k},\frac{C_0^2}{\sigma^2})$
meaning that we would like to bound with $0.05$ (for instance)
$$\pr\{\inf_{0\leq s\leq T}C_s\leq 0\}=\int_0^T\frac{C_0}{\sigma\sqrt{2\pi t^3}}
\exp(-\frac{k^2}{2\sigma^2 t}(t+\frac{C_0}{k})^2)dt=
\int_0^T\frac{C_0}{\sigma\sqrt{2\pi t^3}}
\exp(-\frac{(kt+C_0)^2}{2\sigma^2 t})dt.
$$
After the change of variable $x^2=\frac{C_0^2}{t\sigma^2}$, since $k>0$ this probability is bounded by
$$\int_0^T\frac{C_0}{\sigma\sqrt{2\pi t^3}}
\exp(-\frac{C_0^2}{2\sigma^2 t}))dt=2\int_{C_0/(\sigma\sqrt T)}^\infty\frac{1}{\sqrt{2\pi }}
\exp(-\frac{x^2}{2})dx=\Phi(C_0/(\sigma\sqrt T))$$
where $\Phi$ is the distribution function of the standard Gaussian law.
A sufficient condition to get the cost non-negative on time interval $[0,T]$
with at least probability $0.95$ is
$$C_0/(\sigma\sqrt T)\geq 1.96.$$


\begin{thebibliography}{2001}

\bibitem{Auriol} E. Auriol, P.M. Picard. A theory of BOT concession contracts. {\it Journal of of Economic Behaviour and Organization, 2834}, (2011).

\bibitem {Alain} A. Bensoussan. Stochastic control of partially observable systems,
{\it Cambridge University Press}, (1992).

\bibitem{Biais} B. Biais, T. Mariotti, J.C. Rochet, S. Villeneuve.  Large risks, limited liability and dynamic moral hazard. {\it Econometrica}, Vol. 78, No. 1 (January, 2010), 73--118.

%\bibitem{castaing2} C. Castaing. Sur les multi-applications mesurables {\it Revue fran\c{c}aise
%d'informatique et de recherche opérationnelle}, tome 1, p. 91-126, (1967).



\bibitem{bookCvitanic} J. Cvitanic, J. Zhang.
Contract Theory in Continuous-Time Models, {\it Springer} (2013).

%\bibitem {Aliprantis} Charalambos D. Aliprantis Kim C. Border. Infinite Dimensional Analysis {\it Springer}, third edition.


\bibitem{Touzi} J.  Cvitanic, D. Possama\"i, N. Touzi.
Moral hazard in dynamic risk management, arxiv:1510.07111 (2015).


\bibitem{CMJ} M. Chesney, M. Jeanblanc, M. Yor.   Mathematical Methods for Financial Markets, {\it Springer} (2009).


\bibitem{NEK}  N. El Karoui. Les aspects probabilistes du contr\^ole stochastique, {\it Saint Flour 1979, Volume 876 of  Lectures Notes In Math.}  (1981).


\bibitem{CaroMo2} G.E. Espinosa, C. Hillairet,
B. Jourdain, M. Pontier.  Reducing the debt: Is it optimal to outsource an
investment? {\it (2016). a voir To appear in Mathematics and Financial Economics.}



\bibitem{FR} W. Fleming and R. Rishel.  Deterministic and Stochastic Optimal Control, Springer-Verlag, (1975).


\bibitem{CaroMo} C. Hillairet, M. Pontier.  A Modelisation of Public Private Parternships with failure time,
{\it Laurent Decreusefond and Jamal Najim ed.
Springer Proceedings in Mathematics and Statistics} Vol 22, 91-117, (2012).

%\bibitem{Chewang} Z. Chen and B. Wang Infinite time interval BSDES and and the convergence of g-martingales,
%  J. Austral. Math. Soc. (Series A) 69 (2000), 187-21.
%
%\bibitem{Gilbarg} D. Gilbarg N. S. Trudining. Elliptic Partial Differential Equations of Second Order {\it Springer-Verlag}, year ??
%Berlin Heidelberg New York

\bibitem{IMP} E.  Iossa, D.  Martimort, J.  Pouyet. Partenariats Public-Priv\'e, quelques r\'eflexions. {\it Revue \'economique}, { 59 (3)} (2008).

%\bibitem{IoRe} E.  Iossa,   P.E. Rey.
%Building reputation for contract renewal: Implications for performance dynamics and contract duration,
% {\it Journal of the European Economic Association, European Economic Association},
% vol. 12(3), pages 549-574, 06 (2014).
%
%
%\bibitem{MJ+MY+MC} JEANBLANC M., YOR M., CHESNEY M.: Mathematical Methods for
%Financial Markets. Springer, Berlin, Heidelberg, New York (2009).

%\bibitem{KaS} I. Karatzas, S.   Shreve. Brownian Motion and Stochastic Calculus,
%{\it Springer}, Berlin, Heidelberg, New York (1988).

\bibitem{Krylov}  N. Krylov.   Nonlinear Elliptic and Parabolic Equation of Second Order, D.Reidel, Boston, (1987).

\bibitem{howard} R. Howard.  Dynamic Programming and Markov Processes. MIT Press, Cambridge, (1960).

%\bibitem{Solonnikov} O. A. Ladyzenskaja, V. A. Solonnikov, N. N. Ural'eeva. Linear and Quasilinear Equations of Parabolic Type,... ?

\bibitem{NAKAO} S. Nakao. On the pathwise uniqueness of solutions of one-dimensional stochastic differential equations. Osaka J. Math., 9, 513-518, (1972).

\bibitem{Pages} H. Pag\`es,  D. Possama\"i.   A mathematical treatment of bank monitoring, {\it Finance and Stochastics}, 18(1), 39--73, (2014).

\bibitem{HP} H. Pham.
  Continuous-time stochastic control and optimization with financial applications, {\it  Series Stochastic Modeling and Applied Probability}, vol 61, Springer, (2009).

%\bibitem{protter} P.  Protter.   Stochastic Integration and Differential Equations, Second Edition, Version 2.1 {\it Springer-Verlag}, Heidelberg (2005).

\bibitem{RevYor} D. Revuz and M. Yor. Continuous Martingales and Brownian Motion, {\it Springer}, 2004.


\bibitem{San} Y. Sannikov. A continuous-time version of the principal-agent problem,
{\it Rev. Econ. Studies 75}  957-984, (2008).

%\bibitem{Villeneuve} S. Villeneuve. Optimal Exit under Moral Hazard. {\it cf slides ...}




\bibitem{Strulovici} B. Strulovici, M. Szydlowski. On the smoothness of value functions and the existence of optimal strategies in diffusion models, {\it Journal of Economic Theory}  1016-1055, (2015).
%\bibitem{Rim} Rim Amami. Impulse control problem with switching technology

\end{thebibliography}
\end{document}